\documentclass[conference]{IEEEtran}
\IEEEoverridecommandlockouts
\usepackage{cite}
\usepackage{amsmath,amssymb,amsfonts,amsthm}
\usepackage{algorithmic}
\usepackage{graphicx}
\usepackage{textcomp}
\usepackage{xcolor}
\def\BibTeX{{\rm B\kern-.05em{\sc i\kern-.025em b}\kern-.08em
    T\kern-.1667em\lower.7ex\hbox{E}\kern-.125emX}}
    
\usepackage{cite}
    
\usepackage[colorlinks]{hyperref}
\usepackage[nameinlink]{cleveref}

\newtheorem{proposition}{Proposition}
\newtheorem{corollary}{Corollary}
\newtheorem{lemma}{Lemma}
\newtheorem{remark}{Remark}
\newtheorem{theorem}{Theorem}    

\begin{document}

\title{Missing Mass in Markov Chains
}

\author{\IEEEauthorblockN{Maciej Skorski}
\IEEEauthorblockA{%of\textit{Department...} \\
\textit{University of Luxembourg}\\
maciej.skorski@gmail.com}
}

\maketitle

\begin{abstract}
The problem of missing mass in statistical inference (posed by McAllester and Ortiz, NIPS'02; most recently revisited by Changa and Thangaraj, ISIT'2019) seeks 
to estimate the weight of symbols that have not been sampled yet from a source.

So far all the approaches have been focused on the IID model which, although overly simplistic, is already not straightforward to tackle. The non-trivial part is in handling correlated events and sums of variables with very different scales where classical concentration inequalities do not yield good bounds.

In this paper we develop the research on missing mass further, solving the problem for Markov chains. It turns out that the existing approaches to IID sources are not useful for Markov chains; we reframe the problem as 
studying the tails of hitting times and finding \emph{log-additive approximations} to them. More precisely, we combine the technique of majorization and certain estimates on set hitting times to show how the problem can be eventually reduced back to the IID case. Our contribution are a) new technique to obtain missing mass bounds - we replace traditionally used negative association by majorization which works for a wider class of processes b) first (exponential) concentration bounds for missing mass in Markov chain models c) simplifications of recent results on set hitting times and d) simplified derivation of missing mass estimates for memory-less sources.
\end{abstract}


\section{Introduction}

\subsection{Missing Mass Problem}
The missing mass problem studies the behavior of \emph{unseen} symbols when sampling from a memory-less source. One wants to estimate the probability of elements that have not been visited during $n$-steps ($n$ subsequent samples from a source). Since so far only IID sources have been studied~\cite{mcallester2000convergence,mcallester2003concentration,berend2013concentration,chandra2019concentration,mossel2019impossibility}, it is natural to extend results to process with memory. In this paper we develop such results for Markov chains.

\subsection{Proof Outline for IID}

For the IID case (\cite{DBLP:conf/nips/McAllesterO02, mcallester2003concentration,berend2013concentration, chandra2019concentration} and followed-up works improving bounds) one proves tail bounds as follows

Consider subsequent symbols $X_1,\ldots,X_n$. The fundamental observations is that  the collection
variables $I(X_i=j)$ indexed by tuples of $i,j$ (indicators of whether we hit symbol $j$ at time $i$) are \emph{negatively associated} (referred to as NA).
This follows in two parts; first this is clearly true for any fixed $i$ (indeed then $\sum_j I(X_i=j) = 1$), and then one uses the fact that NA vectors can be augmented~\cite{dubhashi1998balls}.

Let $\tau_j$ be the moment of hitting a symbol $j$ (may be infinite); any symbol that has not been unseen during $n$ steps satisfies $\tau_j > n$. However $\tau_j > n$ is equivalent to $\sum_{i=1}^{n} I(X_i=j) = 0$.
It is also true that block sums of NA variables are NA, therefore $\sum_{i=1}^{n} I(X_i=j)$ indexed by $j$ are NA. Finally threshold transforms preserve NA and thus the set of events $\sum_{i=1}^{n} I(X_i=j) = 0$, equivalent to $\tau_j > n$ is NA.

Therefore the problem reduces to estimating weighted sum of NA boolean variables. This is another non-trivial task where classical inequalities are to week to produce desired results, as they work best with homogenic variables (weights of similar orders of magnitude).

\subsection{Our Result - Markov Chains}

We first extend the problem statement to Markov chains (or, more generally, stationary sources). Consider a Markov chain $(X_i)_i$ over $m$ states $\{1,\ldots,m\}$ with stationary distribution $\pi$ and some initial starting distribution. Let $\tau_j$ be the first time when $j$ is hit. Fix the run length $n$ and consider 
\begin{align}\label{eq:miss_mass_general}
  \textsc{MissingMass} = \sum_{j=1}^{m} \pi(j)\cdot \{\tau_j > n\}
\end{align}
which indeed extends the IID case. Motivated by seeking for possible extensions, we ask the following question
\begin{quote}
    \textbf{Problem}: Do exponentially strong concentration inequalities hold for \textsc{MissingMass} in \Cref{eq:miss_mass_general}?
\end{quote}
To set up expectations correctly, we note that this problem for Markov chains is \emph{much harder than for IID} sources. First, the \emph{NA property fails}; while in the IID case, due to NA (and quite intuitively), not seeing $B$ increases chances for seeing $A$, for Markov chains this may increase or decrease depending on the topology of the chain. For example consider a walk where $A$ is accessed only (or mostly in terms of probability weights) from $B$. 

Second, for memory-less sources the result is based on \emph{exact formulas} for the tails of hitting times which are straightforward to derive. For Markov chains we don't only expect to have accurate formulas but it is very challenging to obtain good lower bounds. Moreover, in our problem we will have to study hitting times of entire sets not individual points. Indeed, because we would need to bound moments of \Cref{eq:miss_mass_general} which yields
mixed moments of $\sum_{j}\{\tau_j > n\}$, we seek for bounds of the form
$\prod\tau_{j\in J}$. Hitting and commute times are generally better understood for points not sets
~\cite{helali2019hitting,levin2017markov}.

We now present our main result. Below $T$ refers to the hitting time of large sets, the quantity studied in recent works~\cite{griffiths2014tight,oliveira2012mixing} (see \Cref{sec:hit_times}).

\begin{theorem}[Majorization by IID Problem]\label{thm:main}
Let the chain $X$ be and $Q_j$ be independent Bernoulli variables with probability $\mathrm{e}^{-c\cdot n\cdot \pi(j) / T}$ for some absolute constant $c$ and $T$ being the maximum hitting time of sets of probability at least $0.5$. %, defined as the maximum time of hitting any set of probability at least 0.5.  
 Then for any set $J$ and any integer $n\geqslant 1$
\begin{align}\label{eq:hits_majorization}
    \Pr\left[\wedge_{j\in J} \{\tau_j > n\}\right] \leqslant \prod_{j\in J}\Pr[Q_j=1]
\end{align}
In particular for any $s>0$ it holds that
\begin{align}\label{eq:conc_majorization}
    \mathbb{E}\exp\left(s\cdot \textsc{MissingMass}\right) \leqslant \mathbb{E}\exp\left(s\cdot \sum_{j=1}^mn \pi(j)Q_j\right).
\end{align}
\end{theorem}

\begin{remark}[Dependency on Hitting Time of Large Sets]
Intuitively the dependency on $T$ is justified, because slow-mixing chains which have large $T$ should have heavy tails for the missing mass.
\end{remark}

\subsection{Discussion and Applications}

\subsubsection{Exponential Bounds for Markovian Sources}

By reducing to the IID case and using bounds from the literature
we obtain (see \Cref{sec:proofs} for a proof)
%~\cite{impagliazzo2010constructive} we obtain the following
\begin{corollary}[Exponential Upper Tails for MCs]\label{cor:exp_bounds_missmass_mcs}
We have the bound $\textsc{MissingMass} \leqslant \mathbb{E}[\sum_j \pi(j) Q_j] + \epsilon$ with probability
at $1-\mathrm{e}^{-\Omega(n\epsilon^2/T)}$.
\end{corollary}
%Note that the bound is one-sided, but as we discussed establishing lower bounds seems very challenging.

\subsubsection{Exponential Bounds for IDD Sources}

Under the IID assumption the expression in \Cref{eq:hits_majorization} can be computed exactly, so that we can actually take $\Pr[Q_j=1]:= \Pr[\tau_j > n] = (1-\pi(j))^{n} \approx \mathrm{e}^{-\pi(j) \cdot n}$. 
This corresponds to setting $c=1$ and $T=1$ in \Cref{thm:main}.
Note that variables on \Cref{eq:hits_majorization} 
Note that sums considered in \Cref{eq:conc_majorization} have equal means by definition of $Q_j$. Thus we re-obtain the same bounds as for the IID case.
\begin{corollary}[Exponential Upper Tails for Missing Mass of IID]
Let $M =   \textsc{MissingMass}$ then under IID $|M - \mathbb{E}M | \leqslant \epsilon$ holds with probability $1-\mathrm{e}^{-\Omega( n\epsilon^2)}$.
\end{corollary}

\subsubsection{Set Hitting Times Estimates}\label{sec:set_hit_estimates}
The bound in \Cref{eq:hits_majorization} is actually the estimate on set hitting times: 
$\wedge_{j\in J} \{\tau_j > n\}$ is equivalent to $\tau_J > n$. Thus we have proven an exponential tail
$\Pr[\tau_J > n]\leqslant \mathrm{e}^{-c\cdot n\cdot \pi(J)/T}$ or, up to a constant, $\mathbb{E}\tau_J = O(T / \pi(J)) $ (these conditions are equivalent up to a constant, see~\Cref{prop:hit_exponential_tail}). See also \Cref{cor:exp_hitting_times}.

\subsubsection{Eliminating Negative Association Theory}
Traditionally the proofs for the IID case depend on non-trivial facts on negative association, for example~\cite{mcallester2003concentration} relies on~\cite{dubhashi1998balls}. However in view of
\Cref{eq:conc_majorization} if the exponential method is used (which is the case of all known bounds) we just need to prove that \Cref{eq:hits_majorization} holds with $\Pr[Q_j] = \Pr[\tau_j > n]$. Plugging this we conclude that one needs to show
\begin{align*}
    \Pr[\tau_J > n] \leqslant \prod_{j\in J} \Pr[\tau_j > n]
\end{align*}
We calculate that $\Pr[\tau_J > n] = (1-\pi(J))^n$ for any $J$, in particular also $\Pr[\tau_j > n] = (1-\pi(j))^n$. Thus it suffices to prove that
\begin{align}
   1-\pi(J) \leqslant \prod_{j\in J} (1-\pi(j))
\end{align}
This follows because $\pi(J) = \sum_{j}\pi(j)$ and from the elementary inequality $(1-a)(1-b)\geqslant 1-a-b$ (applied recursively) valid for all $a,b\in [0,1]$.

\section{Preliminaries}\label{sec:prelim}

We consider a Markov chain $X_1,X_2,\ldots$ over a finite state space $\mathcal{X}$. We assume it is irreducible so that it has a unique stationary distribution $\pi$~\cite{levin2017markov}.

\subsection{Hitting Times}\label{sec:hit_times}
By $T(x,B)$ we denote the expected hitting time of the set $B$ when the chain starts from $x$.
By $T^{+}(A,B)$ we denote the maximal expected hitting time of $B$ over all possible starts in $A$, that is $T^{+}(A,B) = \max_{x\in A} T(A,B)$. Similarly $T^{-}(A,B)$ stands for the minimal hitting time of $B$ over possible starts in $A$, that is $T^{-}(A,B) = \min_{x\in A} T(A,B)$.

We also let $T(B) = T^{+}(\mathcal{X},B)$ (the worse-case expected hitting time of $B$) and
consider the worse expected hitting time to sets of measure at least $\epsilon$, that is 
$T(\epsilon) = \max_{x}\max_{B:\pi(B)\geqslant \epsilon} T(B)$ (here $\pi$ is the stationary distribution).
In our applications we think of $\epsilon$ as a constant and of $T(\epsilon)$ as the hitting time of large sets.

It is a standard fact that for irreducible chains the tails of hitting times are exponential~\cite{roch2015modern,bremaud2017discrete,levin2017markov}. This is shown by splitting long paths into chunks of equal size and applying the Markov property.
\begin{proposition}[Exponential Tails of Hitting Times]\label{prop:hit_exponential_tail}
Fix some initial distribution and let $N_B$ be the hitting time (random variable) of the set $B$. Then we
have $\Pr[N_B > t] \leqslant \exp( -\lceil t / \lfloor \mathrm{e} \cdot \mathbb{E}N_B \rceil \rfloor )$; note that $\mathbb{E} N_B\leqslant T(B)$.
\end{proposition}

\subsection{Ergodicity}

Below we recall the ergodic theorem for Markov Chains~\cite{levin2017markov}
\begin{proposition}[Ergodic Theorem for MCs]\label{thm:ergodic}
If $(X_n)_n$ is an irreducible Markov chain with stationary distribution $\pi$ then
\begin{align}
    \frac{1}{n}\sum_{i=1}^{n} f(X_i) \overset{n\to\infty}{\longrightarrow} \mathbb{E}_{\pi} f\quad\textrm{a.s.}
\end{align}
for any starting distribution of the chain and any real function on the chain states.
\end{proposition}
\subsection{Relative Entropy}
By $D(p\parallel q)$ we note the binary relative-entropy function (Kullback-Leibler divergence), defined as 
$D(p\parallel q) = p\log \frac{p}{q} + (1-p)\log\frac{1-p}{1-q}$. It appears in concentration bounds (aka Chernoff Bounds), and can be bounded from below by the total variation distance as follows~\cite{csiszar2011information}
\begin{proposition}[Pinsker's Inequality]
For any $p,q\in (0,1)$ we have $D(p\parallel q) \geqslant 2(p-q)^2$.
\end{proposition}
 
\section{Proofs}\label{sec:proofs}

\subsection{Bound on Hitting Times}

In this section we prove the bound $\mathbb{E} \tau_J  = O(T/\pi(J))$ which implies then the first part in \Cref{thm:main}, namely
\Cref{eq:hits_majorization}, as discussed in \Cref{sec:set_hit_estimates}.

We will need the following result, which connect the time of reaching $B$ from $A$ to that of the opposite direction, that appears in~\cite{griffiths2014tight}.
The original proof in the arxiv version~\cite{griffiths2012tight_arxiv} was quite involved, based on martingales and concentration inequalities applied in a non-standard setup (the martingale difference were not bounded) and in the final version got subsumed by an argument credited to Peres and Sousi~\cite{griffiths2014tight}.

\begin{lemma}[Bounds on Set Hitting Times]\label{lemma:bound_hit_times}
For an irreducible chain with stationary distribution $\pi$ and any subsets of states $A,B$ we have
\begin{align}\label{eq:bound_hit_times}
    \pi(A) \leqslant \frac{T^{+}(A,B)}{T^{+}(A,B) + T^{-}(B,A)}
\end{align}
In particular
\begin{align}\label{eq:bound_hit_times_2}
 \pi(A)\cdot T^{-}(B,A) \leqslant T^{+}(A,B).
\end{align}
\end{lemma}

Below, as a contribution of independent interest, we provide an alternative simple proof which resembles the approach taken  in~\cite{griffiths2014tight} but uses only simple stopping times rather than martingales and doesn't need concentration inequalities. Before discussing the details we highlight intuitions as follows: we look at how the chain commutes between sets $A$ and $B$. We split the long runs of the chain into rounds where each round is one "return trip": starting from $A$, passing through $B$ and finally returning to $A$. For $m$ rounds on average
we do at least $m\cdot( T^{+}(A,B) + T^{-}(B,A) )$ steps and spent in $A$ at most $m\cdot T^{+}(A,B)$ steps on average. This can be compared to $\pi(A)$ which is the fraction of time spent in $A$ by the Ergodic Theorem (see \Cref{thm:ergodic}).

\begin{proof}
Suppose that the chain starts at some fixed point $x\in A$. >
For $j=1,\ldots,m$ let $N^{A\to B\to A}_{j}$ is the number of further steps it takes the walk to start from $A$ visit $B$ and return to $A$ (such quantities are sometimes called the commute time~\cite{}); also let $N^{A\to B}_{j}$ be the number of further steps it takes the walk to visit $B$ when it starts from $A$.

Once $B$ gets visited the chain will not go into $A$ in this round. Thus we have
\begin{align}
    \textrm{time the walk spent in } A \leqslant \sum_{j=1}^{m} N^{A\to B}_{j}
\end{align}
and clearly
\begin{align}
    \textrm{total time} = \sum_{j=1}^{m}N^{A\to B\to A}_{j}.
\end{align}

To finish the argument we only use convergence in probability. From the discussion it follows that
\begin{align}\label{eq:spent_in_A}
    \pi(A) \leqslant \liminf_{m}{\frac{\sum_{j=1}^{m} N^{A\to B}_{j}}{\sum_{j=1}^{m}N^{A\to B\to A}_{j}}}
\end{align}

In the next step we will replace the numerator and the denominator with their means. For any $\epsilon>0$ and sufficiently big $m$ we have
\begin{align}
    (\pi(A)-\epsilon)\cdot \sum_{j=1}^{m}N^{A\to B\to A}_{j}  \leqslant \sum_{j=1}^{m} N^{A\to B}_{j}.
\end{align}
Note that both collections $\{N^{A\to B}_{j}\}_j$ and $\{N^{A\to B\to A}_{j}\}_j$ are independent as it follows from the Markov property; they are however not identically distributed because of evolving start points. Taking the expectation we obtain
\begin{align}
    (\pi(A)-\epsilon)\cdot \sum_{j=1}^{m}\mathbb{E} N^{A\to B\to A}_{j} \leqslant \sum_{j=1}^{m} \mathbb{E} N^{A\to B}_{j}
\end{align}
Without loosing generality we can assume that $\pi(A) > 0$ and that $\epsilon < \pi(A)$; then $\pi(A)-\epsilon > 0$.We start with $N^{A\to B\to A}_{j}  = N^{A\to B}_j + \left(N^{A\to B\to A}_{j} - N^{A\to B}_j\right)$, that is splitting the commute time at the moment of approaching $B$ and the way back to $A$. Then we have
$\mathbb{E}\left[ N^{A\to B\to A}_{j} - N^{B\to A}_j \right] \geqslant T^{-}(B,A)$. Thus
\begin{align}
    (\pi(A)-\epsilon)\cdot \sum_{j=1}^{m}\left(\mathbb{E} N^{A\to B}_{j} + T^{-}(B,A)\right) \leqslant \sum_{j=1}^{m} \mathbb{E} N^{A\to B}_{j}
\end{align}
Rearranging the terms we write
\begin{align}
    (\pi(A)-\epsilon)\cdot \sum_{j=1}^{m} T^{-}(B,A) \leqslant (1-\pi(A)+\epsilon)\cdot \sum_{j=1}^{m} \mathbb{E} N^{A\to B}_{j}
\end{align}
We now use the bound $\mathbb{E} N^{A\to B}_j \leqslant T^{+}(B,A)$ to arrive at
\begin{align}
    (\pi(A)-\epsilon)\cdot \sum_{j=1}^{m} T^{-}(B,A) \leqslant (1-\pi(A)+\epsilon)\cdot \sum_{j=1}^{m} \mathbb{E} T^{+}(A,B)
\end{align}
This is equivalent to
\begin{align}
    \pi(A) \leqslant \frac{T^{+}(A,B)}{T^{+}(A,B) + T^{-}(B,A)} + \epsilon
\end{align}
and the result follows as $\epsilon$ can be arbitrary small.

We also comment how the basic Chebyszev inequality can be also used to accomplish the argument.
We observe that $N^{A\to B}_{j}$ and $N^{A\to B\to A}_{j}$ have bounded moments, because the tails of stopping times are exponential (see \Cref{prop:hit_exponential_tail}). Then by the Chebyszev inequality implies
\begin{align}
    \Pr\left[ \left|\frac{\sum_{j=1}^{m} \left[ N^{A\to B}_{j}-\mathbb{E}N^{A\to B}_{j} \right]}{m}\right| \geqslant\epsilon \right] = O(1/m\epsilon^2) \\
    \Pr\left[ \left|\frac{\sum_{j=1}^{m} \left[ N^{A\to B\to A}_{j}-\mathbb{E}N^{A\to B \to A}_{j} \right]}{m}\right| \geqslant\epsilon \right] = O(1/m\epsilon^2)
\end{align}
where the constant depends on $A,B$ and the chain. Using this in \Cref{eq:spent_in_A} arrive at the same conclusion.
\end{proof}
\begin{remark}\label{remark:measure_vs_hit_best}
By refining the current proof we can show a slightly better bound with the constant $2$.
\end{remark}

\begin{lemma}[Measure vs Hitting Time]\label{lemma:measure_vs_hit}
For any $A$ we have that $T(A) \leqslant  2 \cdot T(0.5) / \pi(A)$.
\end{lemma}
Before giving a proof we explain the intuition. Consider the set of starting points $B$ that are "unlucky" for $A$, that is make the hitting time very long. Then $T^{-}(B,A)$ is very large and to keep the right-hand side of \Cref{lemma:bound_hit_times}
big enough $T^{+}(A,B)$ must be sufficiently big; more precisely at least by a factor of $1/\pi(A)$. 
But we bounded the hitting time of big sets $B$ (see the definition of $T(\epsilon)$), therefore we conclude that $B$ is small. In other words, the complementary set $B^c$ of good starting points is big and the walk quickly approaches it; and once it gets there it also quickly approaches $A$ by the definition of good starting points.
\begin{proof}
\Cref{eq:bound_hit_times} implies that %$T^{-}(B,A) \leqslant \frac{T^{+}(A,B)}{\pi(A)}$. By definition $T^{+}(A,B) \leqslant T(B)$ also 
$T^{-}(B,A) \leqslant \frac{T(B)}{\pi(A)}$. Let $B$ contain all $x$ such that $T(x,A) > T(0.5)/\pi(A)$ (unlucky starting points). Then we must have $\pi(B)\leqslant 0.5$. Then $\pi(B^c) \geqslant 1-0.5=0.5$ which implies $T(B^c) \leqslant T(0.5)$, 
and $B^c$ are good starts for $A$ that is 
$T^{+}(B^c,A) \leqslant T(0.5)/\pi(A)$. 

By the Markov property
\begin{align*}
T(x,A) \leqslant T(x,B^c) + T^{+}(B^c,A))
\end{align*}
Taking the maximum over $x$ on the right-hand side and using the previous bounds we obtain
\begin{align*}
T(x,A) \leqslant T(0.5) + T(0.5) / \pi(A)
\end{align*}
again taking the maximum over $x$ on the left-hand side we finish the proof.
\end{proof}

Combining \Cref{prop:hit_exponential_tail} and \Cref{lemma:measure_vs_hit} we obtain the following
\begin{corollary}[Explicit Exponential Tails of Hitting Times]\label{cor:exp_hitting_times}
Let $N_A$ be the hitting time of a set $A$ for some initial distribution of the chain. Then we have
\begin{align}
    \Pr[N_A > t] \leqslant \exp( -\Omega(t \cdot \pi(A) / T(0.5)  ))
\end{align}
for some absolute constant under $\Omega(\cdot)$.
\end{corollary}
\begin{remark}[Explicit Constant]
The explicit constant can be set to $1/\mathrm{e}\cdot (1+o(1))$ for large $t$. This follows from \Cref{prop:hit_exponential_tail} and \Cref{remark:measure_vs_hit_best}.
\end{remark}

\subsection{Combining with IID Bounds}

The condition \Cref{eq:hits_majorization}, proved in the previous subsection,
implies \Cref{eq:conc_majorization}. Indeed, if $\prod_{j\in J} u_j \leqslant \prod_{j\in J} q_j$ (majorization) then  $\mathbb{E}\left(\sum_j u_j\right)^k\leqslant \mathbb{E}\left(\sum_j q_j\right)^k$ and, by the Taylor expansion of $\exp(\cdot)$, we obtain $\mathbb{E} \exp(s\cdot \sum_j u_j) \leqslant \mathbb{E} \exp(s\cdot \sum_j u_j)$. This finishes the proof of~\Cref{thm:main}.
Therefore
Therefore upper obtained through the exponential method for IID variables $Q_j$ will apply as well. Following the discussion in~\cite{mcallester2003concentration} (particularly Lemmma 11) we obtain
the upper bound $\sum_j \pi(j)\cdot\mathbb{E}Q_j + \epsilon$
with probability $1-\mathrm{e}^{-\Theta( n\epsilon^2)}$.

\label{eq:exp_bounds_missmass_mcs}

\section{Conclusion}

We have studied the missing mass problem under Markov chain models. The obtanined reduction allows for deriving bounds from an IID scenario.

\section*{Acknowledgment}

%The author thanks to

\bibliographystyle{amsplain}
\bibliography{citations}

\providecommand{\bysame}{\leavevmode\hbox to3em{\hrulefill}\thinspace}
\providecommand{\MR}{\relax\ifhmode\unskip\space\fi MR }
\providecommand{\MRhref}[2]{%
  \href{http://www.ams.org/mathscinet-getitem?mr=#1}{#2}
}
\providecommand{\href}[2]{#2}
\begin{thebibliography}{10}

\bibitem{berend2013concentration}
Daniel Berend, Aryeh Kontorovich, et~al., \emph{On the concentration of the
  missing mass}, Electronic Communications in Probability \textbf{18} (2013).

\bibitem{bremaud2017discrete}
P.~Br{\'e}maud, \emph{Discrete probability models and methods: Probability on
  graphs and trees, markov chains and random fields, entropy and coding},
  Probability Theory and Stochastic Modelling, Springer International
  Publishing, 2017.

\bibitem{chandra2019concentration}
Prafulla Chandra and Andrew Thangaraj, \emph{Concentration and tail bounds for
  missing mass}, 2019 IEEE International Symposium on Information Theory
  (ISIT), IEEE, 2019, pp.~1862--1866.

\bibitem{csiszar2011information}
I.~Csisz{\'a}r and J.~K{\"o}rner, \emph{Information theory: Coding theorems for
  discrete memoryless systems}, Cambridge University Press, 2011.

\bibitem{dubhashi1998balls}
Devdatt Dubhashi and Desh Ranjan, \emph{Balls and bins: A study in negative
  dependence}, Random Structures \& Algorithms \textbf{13} (1998), no.~2,
  99--124.

\bibitem{griffiths2014tight}
Simon Griffiths, Ross Kang, Roberto Oliveira, and Viresh Patel, \emph{Tight
  inequalities among set hitting times in markov chains}, Proceedings of the
  American Mathematical Society \textbf{142} (2014), no.~9, 3285--3298.

\bibitem{griffiths2012tight_arxiv}
Simon Griffiths, Ross~J. Kang, Roberto~Imbuzeiro Oliveira, and Viresh Patel,
  \emph{Tight inequalities among set hitting times in markov chains}, 2012.

\bibitem{helali2019hitting}
Amine Helali and Matthias L{\"o}we, \emph{Hitting times, commute times, and
  cover times for random walks on random hypergraphs}, Statistics \&
  Probability Letters (2019).

\bibitem{levin2017markov}
David~A. Levin, Yuval Peres, and Elizabeth~L. Wilmer, \emph{{Markov chains and
  mixing times}}, American Mathematical Society, 2006.

\bibitem{mcallester2003concentration}
David McAllester and Luis Ortiz, \emph{Concentration inequalities for the
  missing mass and for histogram rule error}, Journal of Machine Learning
  Research \textbf{4} (2003), no.~Oct, 895--911.

\bibitem{DBLP:conf/nips/McAllesterO02}
David~A. McAllester and Luis~E. Ortiz, \emph{Concentration inequalities for the
  missing mass and for histogram rule error}, Advances in Neural Information
  Processing Systems 15 [Neural Information Processing Systems, {NIPS} 2002,
  December 9-14, 2002, Vancouver, British Columbia, Canada], 2002,
  pp.~351--358.

\bibitem{mcallester2000convergence}
David~A McAllester and Robert~E Schapire, \emph{On the convergence rate of
  good-turing estimators.}, COLT, 2000, pp.~1--6.

\bibitem{mossel2019impossibility}
Elchanan Mossel and Mesrob Ohannessian, \emph{On the impossibility of learning
  the missing mass}, Entropy \textbf{21} (2019), no.~1, 28.

\bibitem{oliveira2012mixing}
Roberto Oliveira et~al., \emph{Mixing and hitting times for finite markov
  chains}, Electronic Journal of Probability \textbf{17} (2012).

\bibitem{roch2015modern}
Sebastien Roch, \emph{Modern discrete probability: An essential toolkit},
  Lecture notes (2015).

\end{thebibliography}

\end{document}